\documentclass[11pt, a4paper]{article}
\usepackage{amsmath,amssymb,latexsym,graphics,epsfig}
\usepackage{hyperref}
\usepackage{color}
\usepackage{amsthm}

\newtheorem{teorema}{Theorem}

\newtheorem{proposicio}{Proposition}

\def\prova{{\boldmath  $Proof.$}\hskip 0.3truecm}

\def\final{\mbox{ \quad $\Box$}}

\def\ds{\(\displaystyle}
\def\R{\mathbb R}

\def\exc{\mbox{$\varepsilon$}}
\def\excC{\mbox{$\varepsilon_{\Cpetita}$}}
\def\excCk{\mbox{$\varepsilon_{\Cpetita_k}$}}

\def\e{\mbox{\boldmath $e$}}

\def\vecnu{\mbox{\boldmath $\nu$}}
\def\vecrho{\mbox{\boldmath $\rho$}}

\def\vece{\mbox{\boldmath $e$}}
\def\vecu{\mbox{\boldmath $u$}}

\def\vecz{\mbox{\boldmath $z$}}

\def\vecj{\mbox{\boldmath $j$}}

\def\matrixA{\mbox{\boldmath $A$}}
\def\A{\mbox{\boldmath $A$}}

\def\matrixE{\mbox{\boldmath $E$}}
\def\E{\mbox{\boldmath $E$}}

\def\matrixI{\mbox{\boldmath $I$}}

\def\ev{\mbox{\rm ev}}
\def\sp{\mbox{\rm sp}}
\def\spG{\mbox{\rm sp\,$\Gamma$}}
\def\evG{\mbox{\rm ev\,$\Gamma$}}

\def\spC{\mbox{$\rm sp_{\Cpetita}$\,$\Gamma$}}

\def\evC{\mbox{$\rm ev_{\Cpetita}$\,$\Gamma$}}
\def\evCk{\mbox{$\rm ev_{\Ckpetita}$\,$\Gamma$}}
\def\evCkk{\mbox{$\rm ev_{\Ckkpetita}$\,$\Gamma$}}
\def\evD{\mbox{$\rm ev_{\Dpetita}$\,$\Gamma$}}

\def\Cpetita{\mbox{{\tiny $C$}}}
\def\Ckpetita{\mbox{{\tiny $C_k$}}}
\def\Ckkpetita{\mbox{{\tiny $C_{k+1}$}}}
\def\Dpetita{\mbox{{\tiny $D$}}}
\def\subC{\Cpetita}
\def\subD{\Dpetita}
\def\dsubC{d_{\Cpetita}}

\def\pbar{\overline{p}}

\begin{document}

\title{On the Local Spectra of the Subconstituents  of a Vertex Set and Completely Pseudo-Regular Codes
\footnote{
Research supported by the Spanish Council under project
MTM2011-28800-C02-01 and by the Catalan Research Council under
project 2009SGR01387.}
\author{M. C\'amara, J. F\`abrega, M.A. Fiol, and E. Garriga
\\ \\
{\small Departament de Matem\`atica Aplicada IV} \\
{\small Universitat Polit\`ecnica de Catalunya, BarcelonaTech} \\
{\small Jordi Girona 1-3, 08034 Barcelona, Catalonia} \\
{\small (e-mails: {\tt \{mcamara,jfabrega,fiol,egarriga\}@ma4.upc.edu})}}}

\maketitle

\begin{abstract}
The local spectrum of a vertex  set in a graph has been proven to be
very useful to study some of its metric properties. It also has
applications in the area of pseudo-distance-regularity around a set
and can be used to obtain quasi-spectral characterizations of
completely (pseudo-)regular codes. In this paper we study the
relation between the local spectrum of a vertex set and the local
spectrum of each of its subconstituents. Moreover, we obtain a new
characterization for completely pseudo-regular codes, and
consequently for completely regular codes, in terms of the relation
between the local spectrum of an extremal set of vertices and the
local spectrum of its antipodal set. We also present a new proof of
the version of the Spectral Excess Theorem for extremal sets of
vertices.
\end{abstract}


\section{Introduction}\label{sec: intro}
The notion of local spectrum was first introduced in \cite{fgy96}
for a single vertex of a graph. In that paper, such a concept was
used to obtain several quasi-spectral characterizations of local
(pseudo)-distance-regularity. In the study of
pseudo-distance-regularity around a set of vertices \cite{fg99},
which particularizes to that of completely regular codes when the
graph is regular, the local spectrum is generalized to a set of
vertices. As commented in the same paper, when we study the graph
from a ``base" vertex subset, its local spectrum plays a role
similar to the one played by the (standard) spectrum for studying
the whole graph.

In this work we are interested in the study of the relation between
the local spectrum of a vertex set and the local spectra of the
elements of the distance partition associated to it (also known as
its subconstituents). Thus, Section~\ref{sec: local spectrum} is
devoted to define the local spectrum of a subset of vertices in a
graph. In Section~\ref{sec: codes} completely pseudo-regular codes
are introduced and we discuss some known results of special
interest. Our main results can be found in Sections~\ref{sec: spec
subconstituents} and \ref{sec: charact. codes}, where we give
sufficient conditions implying a tight relation between the local
spectrum of a set of vertices and that of each of its
subconstituents. As a consequence, we obtain a new characterization
of completely (pseudo-)regular codes. In the way we also obtain some
information about the structure of the local spectrum of the
subconstituents associated to a completely pseudo-regular code and
we give a new proof of a result from \cite{fg99}, which can be seen as
the Spectral Excess Theorem \cite{fg97} for sets of vertices.

Before going into our study, let us first give some notation.
In this paper
$\Gamma=(V,E)$ stands for a simple connected graph with vertex set
$V=\{1,2,\ldots,n\}$. Each vertex $i\in V$ is identified with the $i$-th
unit coordinate (column) vector $\e_i$ and ${\cal V}\cong \R^n$ denotes the vector space of formal
linear combinations of its vertices. The adjacencies in $\Gamma$,
$\{i,j\}\in E$, are denoted by $i\sim j$ and $\partial (\cdot,
\cdot)$ stands for the distance function in $\Gamma$. Given a
set of vertices $C\subset V$, the distance from a vertex $i$ to $C$
is given by the expression $\partial(i,C)=\min\{\partial(i,j)\;|\;
j\in C\}$. We denote by $\excC=\max_{i\in V}{\partial(i,C)}$
the eccentricity of $C$. Notice that, in the context of coding
theory, the parameter $\excC$ corresponds to the covering radius of
the code $C$.

As usual, $\A$ stands for the adjacency matrix of $\Gamma$, with set
of different eigenvalues
$\evG :=\ev \A=\{\lambda_0,\lambda_1,\ldots,\lambda_d\}$, where
$\lambda_0>\lambda_1>\cdots>\lambda_d$. The spectrum of $\Gamma$ is
$$
\spG:=\sp \matrixA=\{\lambda_0^{m(\lambda_0)},
\lambda_1^{m(\lambda_1)}, \ldots, \lambda_d^{m(\lambda_d)}\},
$$
where $m(\lambda_l)$ is the multiplicity of the eigenvalue
$\lambda_l$. We denote by ${\cal E}_l=\ker(\A-\lambda_l\matrixI)$
the eigenspace of $\A$ corresponding to $\lambda_l$. Recall that,
since $\Gamma$ is connected, ${\cal E}_0$ is one-dimensional, and all its elements are eigenvectors
having all its components either positive or negative (see
e.g. \cite{biggs,CvDoSa79}). Denote by
$\vecnu=(\nu_1,\nu_2,\ldots, \nu_n)\in{\cal E}_0$ the unique
positive eigenvector of $\Gamma$ with minimum component equal to 1.

Note that ${\cal V}$ is a module over the quotient ring $\R[x]/{(Z)}$, where $(Z)$ is the ideal generated by the minimal polynomial
of $\A$, $Z=\prod_{l=0}^d(x-\lambda_l)$, with product defined by
$$
p\vecu:=p(\matrixA)\vecu\qquad\mbox{\rm for every $p\in \R[x]/{(Z)}$ and
 $\vecu\in {\cal V}$.} $$
With this notation, let us remark that the orthogonal projection
$\matrixE_l$ of ${\cal V}$ onto the eigenspace ${\cal E}_l$ corresponds to
$$
\matrixE_l\vecu=Z_l\vecu,\quad\vecu\in{\cal V},
$$
where $Z_l$, $l=0,1,\ldots,d$, is the Lagrange interpolating polynomial
satisfying $Z_l(\lambda_h)=\delta_{lh}$, that is:
$$
Z_l=\frac{(-1)^l}{\pi_l}\prod_{0\leq h\leq d,\,h\neq l}
(x-\lambda_h),
$$
with $\pi_l$ being the moment-like parameter given by
$\pi_l:=\prod_{0\leq h\leq d,\,h\neq l}|\lambda_l-\lambda_h|$.

\section{$C$-Local spectrum}\label{sec: local spectrum}

Given a set of vertices $C$ of $\Gamma$, define the map
$\vecrho:{\cal P}(V)\mapsto {\cal V}$ by
$\vecrho\emptyset:=\mathbf{0}$ and $\vecrho C:=\sum_{i\in
C}\nu_i\vece_i$ for $C\neq\emptyset$. Consider the spectral
decomposition of the unit vector $\vece_{\Cpetita}=\vecrho C/\|\vecrho
C\|=\vecz_{\Cpetita}(\lambda_0)+\vecz_{\Cpetita}(\lambda_1)
+\cdots+\vecz_{\Cpetita}(\lambda_d)$, that is
$\vecz_{\Cpetita}(\lambda_l)=\matrixE_l\vece_{\Cpetita}\in {\cal
E}_l$, $0\le l\le d$. The $C$-{\it multiplicity} of the eigenvalue
$\lambda_l$ is defined by
$$
m_{\Cpetita}(\lambda_l):=\langle \E_l\vece_{\Cpetita},\vece_{\Cpetita}\rangle = \|\vecz_{\Cpetita}(\lambda_l)\|^2.
$$

If $\mu_0>\mu_1>\cdots >\mu_{\dsubC}$ are the eigenvalues of
$\Gamma$ with nonzero $C$-multiplicity, the $C$-{\it local spectrum}
of $\Gamma$ is defined by
$$
\spC:=\{\mu_0^{m_{\subC}(\mu_0)},\mu_1^{m_{\subC}(\mu_1)},\ldots,
         \mu_{\dsubC}^{m_{\subC}(\mu_{\dsubC})}\},
$$
and we denote by $\evC:=\{\mu_0,\mu_1,\ldots,\mu_{\dsubC}\}$,
$\mu_0>\mu_1>\cdots >\mu_{\dsubC}$, the set of different eigenvalues
in the $C$-local spectrum. Let us remark that, as
$\E_0\vece_{\Cpetita}=\frac{\langle \vece_{\Cpetita},
\vecnu\rangle}{\|\vecnu\|^2}\vecnu=\frac{\|\vecrho
C\|}{\|\vecnu\|^2}\vecnu$, we have
$m_{\Cpetita}(\lambda_0)=\frac{\|\vecrho C\|^2}{\|\vecnu\|^2}\neq
0$, and hence $\mu_0=\lambda_0$. The parameter $\dsubC$ is called
the {\it dual degree} of $C$ and it provides an upper bound for the
eccentricity of the vertex set, $\excC\leq \dsubC$ (see \cite{fg99}). When
the equality is attained we say that $C$ is {\it extremal}.

Consider the idempotents $\matrixE_l^{\Cpetita}$, $0\le l\le
\dsubC$, corresponding to the members of the $C$-spectrum, that is
$\matrixE_l^{\Cpetita}$ is the projection of ${\cal
V}_{\subC}=\bigoplus_{\lambda_h\in\evC}{\cal E}_h$ onto the
eigenspace corresponding to $\mu_l$. As we have done for the
standard spectrum of the graph, we define for each $\mu_l\in\evC$
the moment-like parameter $\pi_l(C):=\prod_{0\leq h\leq \dsubC,\,
h\neq l}|\mu_l-\mu_h|$ and consider the polynomial
\begin{equation}\label{eq polinomis proj C}
Z_l^{\Cpetita}:=\frac{(-1)^l}{\pi_l(C)}\prod_{0\leq h\leq
\dsubC,\,h\neq l}(x-\mu_h)
\end{equation}
which gives $Z_l^{\Cpetita}(\matrixA)=\matrixE_l^{\Cpetita}$.

\section{Completely pseudo-regular codes}\label{sec: codes}
In this section we review some known results on completely
pseudo-regular codes. These results are formulated in terms of
$C$-local pseudo-distance-regularity, which extends the notion of
local distance-regularity from single vertices to subsets of
vertices and from regular to non-necessarily regular graphs.
Completely pseudo-regular codes where introduced in \cite{fg99} with
the aim of generalizing both local distance-regularity and completely
regular codes \cite{g93}. Let us consider the distance partition of a vertex set $C\subset V$ of
$\Gamma$, given by the sets $C_k=\{i\in V \;|\;
\partial(i,C)=k \}$, $k=0,1,\ldots,\excC$, which are known as the
{\em subconstituents} associated to $C$. Consider also the functions
$a,b,c: V\longrightarrow [0,\lambda_0]$ acting on a vertex $i\in
C_k$ as follows:

\[ \begin{array}{l}
     c(i)=\left\{\begin{array}{ll}
         \mbox{\ds 0 \)}& \quad\quad(k=0);\\
         \mbox{\ds
         \frac{1}{\nu_i}\sum_{j\in\Gamma(i)\cap C_{k-1}}\hspace{-5mm}\nu_j\)
              } &
                     \quad\quad(1\leq k\leq \excC).\end{array}\right.\\
         \\
      a(i)=\mbox{\ds \frac{1}{\nu_i}\sum_{j\in\Gamma(i)\cap C_{k}}\)}
                      \hspace{-3mm}\nu_j \quad\quad\quad\quad (0\leq k\leq \excC).  \\
                    \\
      b(i)=\left\{\begin{array}{ll}
             \mbox{\ds \frac{1}{\nu_i}\sum_{j\in\Gamma(i)\cap C_{k+1}}\hspace{-5mm}\nu_j\)}  &
                     \;\quad\quad(0\leq k\leq \excC-1);\\
             \mbox{\ds 0 \)} & \;\quad\quad(k=\excC). \end{array}\right.
   \end{array}\]
We say that $C$ is a {\it completely pseudo-regular code} (or that $\Gamma$ is {\em $C$-local pseudo-distance-regular}) when the values of $c$, $a$ and $b$ do not depend on the chosen vertex $i\in C_k$ but only on $k$; that is, when the distance
partition with respect to $C$ is pseudo-regular (see \cite{f99}). If
this is the case, we denote by $c_k$, $a_k$ and $b_k$ the values
of these three functions on the vertices of $C_k$,
$k=0,1,\ldots,\excC$, and refer to them as the $C$-{\it local pseudo-intersection numbers}.

Remark that, since $\vecnu\in{\cal E}_0$, the sum of these three
functions is constant over all the vertices of $\Gamma$:
$$
c(i)+a(i)+b(i)=\frac{1}{\nu_i}\sum_{j\in\Gamma(i)}\nu_j=\frac{1}{\nu_i}(\A\vecnu)_i=\lambda_0 \qquad \mbox{for all } i\in V.
$$
Thus, in a sense, we can say that the weight $\nu_i$, given to every
vertex $i$, regularizes the graph. Note also that, when $\Gamma$ is a
regular graph we have that $\vecnu=\vecj$, the all-$1$ vector, and the
definition of a completely pseudo-regular code particularizes to
that of a completely regular code.

\subsection{Some characterizations of $C$-local
pseudo-distance-regularity}\label{subsec: codes charact}

As established in \cite{cffg10}, one can consider two approaches to
completely pseudo-regular codes. One is based on the
combinatorial definition given in the previous section and the other
relies upon the study of the Terwilliger algebra \cite{t92} associated
to a vertex set. Although the results presented here can be
obtained from both points of view, we will focus the attention on
the combinatorial approach, which makes more evident the key role of the
local spectrum.

From now on, all the polynomials must be considered in the
quotient ring $\R[x]/(Z_{\Cpetita})$, where
$Z_{\Cpetita}=\prod_{\mu_l\in \evC}(x-\mu_l)$. We begin by defining
the $C$-{\it local scalar product} as follows:
$$
\langle p,q\rangle_{\Cpetita}:=\langle p\vece_{\Cpetita},
q\vece_{\Cpetita}\rangle=\sum_{l=0}^{d_{\Cpetita}}
m_{\Cpetita}(\mu_l)p(\mu_l)q(\mu_l).
$$
Then, we consider the family of {\em $C$-local predistance polynomials}, $(p_k)_{0\leq k\leq
d_{\Cpetita}}$, which is the unique orthogonal system with respect
to the $C$-local scalar product, such that $\deg p_k =k$ and
$\langle p_k,p_h\rangle_{\subC}=\delta_{kh}p_k(\mu_0)$,
$0\le k\le d_{\subC}$; see \cite{cffg09}. In \cite{fg99} one can
find some quasi-spectral characterizations of $C$-local
pseudo-distance-regularity, which are based on these
polynomials. In particular, the following two theorems are of special interest
in our work.

\begin{teorema}
\label{teo1 pdrC} Let $\Gamma=(V,E)$ be a graph and let $C\subset V$ have
eccentricity $\excC$. Then, $\Gamma$ is $C$-local
pseudo-distance-regular if and only if the $C$-local
$($pre$)$distance polynomials satisfy $\vecrho C_k=p_k\vecrho C$,
$k=0,1\ldots,\excC$.
\end{teorema}

The following result establishes that, if we assume that $C$ is
extremal, it is enough to check the condition of Theorem \ref{teo1
pdrC} for the set of vertices at maximum distance $C_{\excC}$. We
refer to this set as the {\it antipodal} set of $C$ and, if
there is no possible confusion, we  write $D=C_{\excC}$.

\begin{teorema}\label{teo2 pdrC}
Let $\Gamma=(V,E)$ be a graph and let $C$ be an extremal set with eccentricity $\excC=\dsubC$. Then, $\Gamma$
is $C$-local pseudo-distance-regular if and only if there exists a
polynomial $p\in \R_{\dsubC}[x]$ such that $p\vecrho C=\vecrho D$, in
which case $p$ is the $C$-local predistance
polynomial $p_{d_{\Cpetita}}$.
\end{teorema}

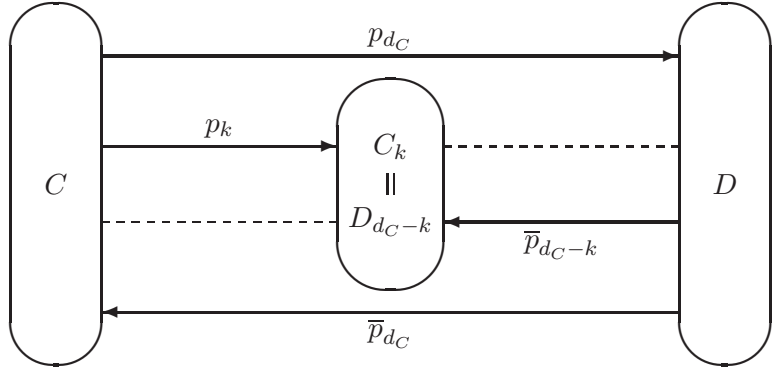
\begin{figure}
\begin{center}
\unitlength=1mm
\begin{picture}(100,50)(-50,-25)
\put(0,0){{
\thicklines{ \put(-44,0){\oval(12,48)} \put(0,0){\oval(14,28)}
\put(44,0){\oval(12,48)} \put(-38,17){\vector(1,0){76}}
\put(38,-17){\vector(-1,0){76}} \put(-38,5){\vector(1,0){31}}
\put(38,-5){\vector(-1,0){31}} } \put(-44,0){\makebox(0,0){$C$}}
\put(44,0){\makebox(0,0){$D$}} \put(0,5){\makebox(0,0){$C_k$}}
\put(0,-5){\makebox(0,0){$D_{\dsubC-k}$}}
\put(-22.5,7.4){\makebox(0,0){$p_k$}}
\put(22.5,-8){\makebox(0,0){$\pbar_{\dsubC-k}$}}
\put(0,19.2){\makebox(0,0){$p_{\dsubC}$}}
\put(0,-19.8){\makebox(0,0){$\pbar_{\dsubC}$}} \thinlines{
\put(-0.4,1.4){\line(0,-1){2.8}} \put(0.4,1.4){\line(0,-1){2.8}}
\put(7,5){\dashbox{1}(31,0)} \put(-38,-5){\dashbox{1}(31,0)} } }}
\end{picture}
\end{center}
\caption{\label{fig-tight} Antipodal completely pseudo-regular codes
and their distance partition}
\end{figure}

This last result points out the relevance of the antipodal set in
the study of completely pseudo-regular codes. In fact, it is known
from \cite{fg99} that a set of vertices is a completely
pseudo-regular code if and only if its antipodal set is.
This symmetry is illustrated in Fig.~\ref{fig-tight} where the $\pbar_{k}$'s are the $D$-local predistance polynomials.

\section{Spectra of the subconstituents}\label{sec: spec subconstituents}
As showed in \cite{cffg10}, for an extremal set of vertices $C$ with
antipodal set $D$, the relation between the local spectra of $C$
and $D$ is tight and can be expressed in terms of the following
proposition. We include the proof for the sake of completeness.
\begin{proposicio} \label{prop desigualtat producte multiplicitats}
Let $C$ be an extremal set and let $D$ be its antipodal set. Then,
$\evC\subset \evD$ and the $C$-multiplicities and $D$-multiplicities
satisfy
\begin{equation}\label{eq. c,d-multi}
m_{\subC}(\mu_l)m_{\subD}(\mu_l)\geq \frac{\pi_0^2(C)}{\pi_l^2(C)}\,
\frac{\|\vecrho C\|^2\|\vecrho D\|^2}{\|\vecnu\|^4}
     \quad\mbox{ for all }\;\mu_l\in\evC.
\end{equation}
Moreover:
\begin{itemize}
\item[$(i)$] Equality holds in Eq. $(\ref{eq. c,d-multi})$ for an
eigenvalue $\mu_l$ if and only if the vectors $\vecz_{\Cpetita}(\mu_l)$ and $\vecz_{\Dpetita}(\mu_l)$ are linear
dependent.
\item[$(ii)$] Equality in Eq. $(\ref{eq. c,d-multi})$ for every $\mu_l\in\evC$ is equivalent to the existence of a polynomial $p\in\R_{\excC}[x]$ such that
$$
\vecrho D=p\vecrho C+\vecz,\quad where \quad
\textstyle \vecz\in \bigoplus_{\lambda_l\in\evD\setminus\evC}{\cal E}_l.
$$
\end{itemize}
\end{proposicio}

\prova
Consider the $C$-local Hoffman polynomial (see \cite{cffg10,h63})
$$
H_{\Cpetita}=\frac{\|\vecnu\|^2}{\|\vecrho
C\|^2}Z_0^{\subC}=\frac{\|\vecnu\|^2}{\pi_0(C)\|\vecrho
C\|^2}\prod_{l=1}^{\dsubC}(x-\mu_l),
$$
which satisfy $H_{\Cpetita}\vecrho C=\vecnu$. Since both
$Z_l^{\subC}$, defined in Eq. (\ref{eq polinomis proj C}), and
$H_{\Cpetita}$ have degree $d_{\Cpetita}$ and their leading
coefficients are, respectively, $\frac{(-1)^l}{\pi_l(C)}$ and
$\frac{\|\vecnu\|^2}{\pi_0(C)\|\vecrho C\|^2}$, the polynomial
$$
T=\pi_0(C)\frac{\|\vecrho
C\|^2}{\|\vecnu\|^2}H_{\subC}-(-1)^l\pi_l(C)Z_l^{\subC}
$$
has degree less than $d_{\Cpetita}=\excC$.
Then, the vectors
$T\vece_{\Cpetita}$ and $\vece_{\Dpetita}$ are orthogonal, giving:
\begin{eqnarray}
          \langle T\vece_{\Cpetita},\vece_{\Dpetita}\rangle & = &
          \pi_0(C)\frac{\|\vecrho C\|^2}{\|\vecnu\|^2}\langle H_{\subC}\vece_{\Cpetita},\vece_{\Dpetita}\rangle-
          (-1)^l\pi_l(C)\langle Z_l^{\subC}\vece_{\Cpetita},\vece_{\Dpetita}\rangle \nonumber\\
          & = &
          \pi_0(C)\frac{\|\vecrho C\|}{\|\vecrho D\|\,\|\vecnu\|^2}\langle H_{\subC}\vecrho C,\vecrho D\rangle-
          (-1)^l\pi_l(C)\langle \vecz_{\subC}(\mu_l),\vece_{\Dpetita}\rangle \nonumber\\
          & = &
         \pi_0(C)\frac{\|\vecrho C\|\,\|\vecrho D\|}{\|\vecnu\|^2}-
          (-1)^l\pi_l(C)\langle \vecz_{\subC}(\mu_l),\vecz_{\subD}(\mu_l)\rangle=0 \label{eq producte escalar de les projeccions},
\end{eqnarray}
since
$\langle H_{\subC}\vecrho C,\vecrho D\rangle = \langle \vecnu, \vecrho D \rangle=\|\vecrho D\|^2$.
Therefore, by the Cauchy-Schwarz inequality,
$$
\frac{\pi_0^2(C)}{\pi_l^2(C)}\,\frac{\|\vecrho
  C\|^2\|\vecrho D\|^2}{\|\vecnu\|^4}=\langle \vecz_{\subC}(\mu_l),\vecz_{\subD}(\mu_l)\rangle^2
  \leq
  m_{\subC}(\mu_l)m_{\subD}(\mu_l),
$$
where equality occurs if and only if $\vecz_{\subC}(\mu_l)$ and
$\vecz_{\subD}(\mu_l)$ are colinear, and $(i)$ is also proved.
Moreover, as all the terms involved are positive,
$m_{\subD}(\mu_l)>0$ and hence $\mu_l\in \evD$.

In order to proof $(ii)$, suppose now that the equality holds for
every eigenvalue of the $C$-local spectrum. Then, given $\mu_l\in
\evC$, the vectors $\vecz_{\subD}(\mu_l)$, $\vecz_{\subC}(\mu_l)$
are proportional. More precisely, by Eq. (\ref{eq producte escalar de
les projeccions}), there exist $\xi_l>0$ such that
$\vecz_{\subD}(\mu_l)=(-1)^l\xi_l\vecz_{\subC}(\mu_l)$. Let $p$ be
the unique polynomial in $\R_{\excC}[x]$ such that
$p(\mu_l)=(-1)^l\frac{\|\vecrho D\|}{\|\vecrho C\|}\xi_l$ for all
$\mu_l\in \evC$. We have
\begin{eqnarray*}
 \matrixE_l\vecrho D&=&\|\vecrho D\|\vecz_{\subD}(\mu_l)=(-1)^l\|\vecrho D\|\xi_l\vecz_{\subC}(\mu_l)\\
 &=&(-1)^l\frac{\|\vecrho D\|}{\|\vecrho C\|}\xi_l\matrixE_l\vecrho C=p(\mu_l)\matrixE_l\vecrho C
 =\matrixE_lp\vecrho C.
\end{eqnarray*}
Thus $\vecz=\vecrho D-p\vecrho
C\in\bigoplus_{\lambda_l\in\evD\setminus\evC}{\cal E}_l$.
Conversely, assuming the existence of $p\in\R_{\excC}[x]$
satisfying $\vecrho D=p\vecrho C+\vecz$,  with $\vecz\in
\bigoplus_{\lambda_l\in\evD\setminus\evC}{\cal E}_l$, by projecting
onto the eigenspace of $\mu_l$ $(\mu_l\in \evC)$ we obtain
$\|\vecrho D\|\vecz_{\subD}(\mu_l)=p(\mu_l)\|\vecrho
C\|\vecz_{\subC}(\mu_l)$ and, by $(i)$, equality in Eq. (\ref{eq.
c,d-multi}) holds for every $\mu_l\in \evC$. \final

Although the last result cannot be directly extended to the
other subconstituents, it suggests that the existence of the distance
polynomials guarantees a tight relation between the local
spectra of the subconstituents. The following proposition
supports this claim.

\begin{proposicio}\label{prop espectresck}
Let $C$ be a completely pseudo-regular code in a graph $\Gamma$. Denote by
$C_k$, $k=0,1,\ldots,\exc(=d_{\subC})$ the subconstituents associated to $C$.
Then $\evCk\subset\evC$. Moreover, for each $l=0,1,\ldots,\dsubC$,
the projections of $\vecrho C_k$ and $\vecrho C$ onto the eigenspace
${\cal E}_l$ are linearly dependent.
\end{proposicio}

\prova
 Let $(p_k)_{0\leq k\leq \exc}$ be the $C$-local predistance
polynomials. From Theorem~\ref{teo1 pdrC} we have that $\vecrho
C_k=p_k\vecrho C$, $k=0,1,\ldots, \exc$.
By projecting onto the eigenspace ${\cal E}_l$ we obtain
$
\matrixE_l\vecrho C_k = \matrixE_l p_k\vecrho C=p_k(\lambda_l)\E_l
\vecrho C
$, so that
$\matrixE_l\vecrho C_k$ and $\matrixE_l\vecrho C$ are colinear. Moreover, since $\vece_{\Cpetita_k}=\frac{\vecrho
C_k}{\|\vecrho C_k\|}$ and $\vece_{\Cpetita}=\frac{\vecrho
C}{\|\vecrho C\|}$  we obtain:
\begin{eqnarray}
m_{\Cpetita_k}(\lambda_l)&=&\|\matrixE_l \vece_{\Cpetita_k}
\|^2=\frac{\|\matrixE_l \vecrho C_k \|^2} {\|\vecrho
C_k\|^2}=\frac{\|p_k(\lambda_l)\matrixE_l\vecrho C\|^2}{\|\vecrho
C_k\|^2}\nonumber\\
&=&\frac{\|\vecrho C\|^2}{\|\vecrho C_k\|^2}(p_k(\lambda_l))^2\|\E_l\vece_{\Cpetita}\|^2\nonumber\\
&=&\frac{\|\vecrho C\|^2}{\|\vecrho C_k\|^2}(p_k(\lambda_l))^2
m_{\Cpetita}(\lambda_l).\label{eq igualtat multiplicitats}
\end{eqnarray}
Thus, $\lambda_l\in \evCk$ and $\evCk\subset\evC$.
\final

Note that,  since for an extremal set we have $\evC\subset\evD$, the
last result shows that, in particular, for a completely
pseudo-regular code we have $\evC=\evD$.

As a by-product, from Eq. (\ref{eq igualtat multiplicitats}),
case $k=\dsubC$, and Proposition~\ref{prop desigualtat producte
multiplicitats}$(ii)$ we have that, for a completely pseudo-regular
code,
\begin{eqnarray}
m_{\Cpetita}(\mu_l)&=&\frac{\pi_0(C)}{\pi_l(C)}\frac{\|\vecrho
D\|^2}{\|\vecnu\|^2}\frac{1}{p_{\dsubC}(\mu_l)};\label{eq md 0}\\
m_{\Dpetita}(\mu_l)&=&\frac{\pi_0(C)}{\pi_l(C)}\frac{\|\vecrho
C\|^2}{\|\vecnu\|^2}p_{\dsubC}(\mu_l)\label{eq md 1};
\end{eqnarray}
for all $\mu_l\in \evC$.  As commented, recall that in this case $D$
is also a completely pseudo-regular code and $C$ is its antipodal
set. If we denote by $\pbar$ the $D$-local predistance polynomial
with maximum degree and use the existing symmetry (see Fig.
\ref{fig-tight}) Eq. (\ref{eq md 0}) gives:
$$
m_{\Dpetita}(\mu_l)=\frac{\pi_0(C)}{\pi_l(C)}\frac{\|\vecrho
C\|^2}{\|\vecnu\|^2}\frac{1}{\pbar(\mu_l)} \qquad \mbox{for all
}\mu_l\in\evC.
$$
This jointly with Eq. (\ref{eq md 1}) leads us to the existing relation
between the $C$-local and $D$-local predistance polynomials of
maximum degree:
$$
\pbar(\mu_l)=\frac{1}{p_{\dsubC}(\mu_l)}\qquad \mbox{ for all
$\mu_l\in\evC$}.
$$

Recall that the $C$-local predistance polynomials satisfy a three
term recurrence (see \cite{cffg09}) and, in particular, for a
completely pseudo-regular code with intersection numbers $a_k$,
$b_k$ and $c_k$, $k=0,1,\ldots,\excC(=\dsubC)$, we have:
$$
xp_k = b_{k-1}p_{k-1}+a_kp_k+ c_{k+1}p_{k+1}\qquad (0\leq k\leq
\dsubC),
$$
where $b_{-1}=c_{\dsubC+1}=0$ and $b_k c_{k+1}>0$ (see \cite{cffg09}).
Thus, for an eigenvalue $\mu_l\in \evC$
$$
p_{k+1}(\mu_l)=\frac{(\mu_l-a_k)p_k(\mu_l)-b_{k-1}p_{k-1}(\mu_l)}{c_{k+1}},
$$

and, since $\evC=\evD$, Eq. (\ref{eq igualtat multiplicitats}) ensures
us that in a completely pseudo-regular code there cannot exist $k$
such that $\mu_l\notin \evCk\cup\evCkk$. That is, for every
$k=0,1,\ldots, \excC-1$, $\evC=\evCk\cup\evCkk$. Remark also that in
this case we have that $\varepsilon_{\Cpetita_1}\geq \dsubC-1$
and, in general, $\excCk\geq\max\{k,\dsubC-k\}$. Thus, since
$d_{\Cpetita_k}\geq\excCk$, the dual degree of the $k$-th
subconstituent satisfies the bound
$$
d_{\Cpetita_k}\geq\max\{k,\dsubC-k\}.
$$
Consequently, in a completely pseudo-regular code, the $C$-local
predistance polynomial $p_k$ vanishes at most at $\min\{k,
\dsubC-k\}$ different eigenvalues of the $C$-local spectrum.

\section{Characterizations of completely pseudo-regular codes}\label{sec: charact. codes}

Some of the results given in the previous section can be used to
obtain characterizations for completely pseudo-regular codes. In
particular, we get an alternative proof for extremal
sets of a result in \cite{fg99}, which can be seen as the version of
the Spectral Excess Theorem \cite{fg97} for sets of vertices.

\begin{teorema} Let $C$ be an extremal set and $D$ its antipodal set. Then
\begin{equation}\label{eq excess}
\frac{\|\vecrho D\|^2}{\|\vecrho
C\|^2}\leq\frac{1/m_{\subC}(\lambda_0)^2
\pi_0^2(C)}{\sum_{l=0}^{\dsubC}1/m_{\subC}(\mu_l) \pi_l^2(C)},
\end{equation}
and the equality holds if and only if $C$ is a completely
pseudo-regular code.
\end{teorema}

\prova By Proposition~\ref{prop desigualtat producte multiplicitats}
we have:
$$
m_{\subD}(\mu_l)\geq \frac{\pi_0^2(C)}{\pi_l^2(C)}\, \frac{\|\vecrho
C\|^2\|\vecrho D\|^2}{\|\vecnu\|^4}\frac{1}{m_{\subC}(\mu_l)}
     \quad\mbox{ for all }\;\mu_l\in\evC.
$$
By adding up for all $\mu_l\in\evC$ and using that we obtain
$m_{\subC}(\lambda_0)=\frac{\|\vecrho C\|^2}{\|\vecnu\|^2}$

\begin{equation}\label{eq sum multiplicitats D}
\sum_{l=0}^{\dsubC}\frac{\pi_0^2(C)}{\pi_l^2(C)}\, \frac{\|\vecrho
D\|^2}{\|\vecrho
C\|^2}\frac{m_{\subC}(\lambda_0)^2}{m_{\subC}(\mu_l)}\leq
\sum_{l=0}^{\dsubC}m_{\subD}(\mu_l)\leq
\sum_{l=0}^{d}m_{\subD}(\lambda_l)=\|\e_{\subD}\|^2=1,
\end{equation}

giving Eq. (\ref{eq excess}). In case of equality in Eq. (\ref{eq sum
multiplicitats D}) we obtain $m_{\subD}(\mu_l)=0$ if
$\lambda_l\notin \evC$, so that $\evC=\evD$. Moreover, in this case
Proposition~\ref{prop desigualtat producte multiplicitats}$(ii)$
applies and there exist a polynomial $p\in\R_{\excC}[x]$ such that
$p\vecrho C=\vecrho D$, or, equivalently, $C$ is a completely
pseudo-regular code. \final

Next theorem gives a new approach to completely pseudo-regular codes
and shows that, when the considered set of vertices is extremal, the
converse of Proposition~\ref{prop espectresck} also holds.

\begin{teorema}
Let $C$ be an extremal set with eccentricity $\exc=d_{\subC}$ and
antipodal set $D=C_{\varepsilon}$. Then $C$ is a completely
pseudo-regular code if and only if the orthogonal projections of
$\vecrho C$ and $\vecrho D$ onto each eigenspace of $\Gamma$ are
colinear.
\end{teorema}

\prova Proposition~\ref{prop espectresck} guaranties the necessity.
Assume now that $\matrixE_l \vecrho D$ and $\matrixE_l \vecrho C$
are colinear for each $l=0,1,\ldots, d$. That is, there exist
constants $\alpha_{l}$ satisfying $\matrixE_l \vecrho
D=\alpha_{l}\matrixE_l\vecrho C$. Let $p_{\exc}$ be the unique
polynomial of degree $\exc$ such that
$p_{\exc}(\lambda_l)=\alpha_{l}$. Then,
$$
\vecrho D=\sum_{l=0}^{d}\matrixE_l \vecrho D = \sum_{l=0}^d
\alpha_{l}\matrixE_l\vecrho C=\sum_{l=0}^{d}p_{\exc}(\lambda_l)
\matrixE_l\vecrho C=p_{\exc}\vecrho C,
$$
and the result follows from Theorem \ref{teo2 pdrC}.
\final

In particular, if the underlying graph $\Gamma$ is regular, the
theorem gives a new characterization of completely regular codes.

\end{document}